
\font\twelverm=cmr12
\font\twelvei=cmmi12
\font\twelvesy=cmsy10
\font\twelvebf=cmbx12
\font\twelvett=cmtt12
\font\twelveit=cmti12
\font\twelvesl=cmsl12
\font\ninerm=cmr9
\font\ninei=cmmi9
\font\ninesy=cmsy9
\font\ninebf=cmbx9
\font\ninett=cmtt9
\font\nineit=cmti9
\font\ninesl=cmsl9
\font\eightrm=cmr8
\font\eighti=cmmi8
\font\eightsy=cmsy8
\font\eightbf=cmbx8
\font\eighttt=cmtt8
\font\eightit=cmti8
\font\eightsl=cmsl8
\font\sixrm=cmr6
\font\sixi=cmmi6
\font\sixsy=cmsy6
\font\sixbf=cmbx6
\font\caps=cmcsc10
\catcode`@=11 
\newskip\ttglue
\def\twelvepoint{\def\rm{\fam0\twelverm}
\textfont0=\twelverm  \scriptfont0=\ninerm  
\scriptscriptfont0=\sevenrm
\textfont1=\twelvei  \scriptfont1=\ninei  \scriptscriptfont1=\seveni
\textfont2=\twelvesy  \scriptfont2=\ninesy  
\scriptscriptfont2=\sevensy
\textfont3=\tenex  \scriptfont3=\tenex  \scriptscriptfont3=\tenex
\textfont\itfam=\twelveit  \def\it{\fam\itfam\twelveit}%
\textfont\slfam=\twelvesl  \def\sl{\fam\slfam\twelvesl}%
\textfont\ttfam=\twelvett  \def\tt{\fam\ttfam\twelvett}%
\textfont\bffam=\twelvebf  \scriptfont\bffam=\ninebf
\scriptscriptfont\bffam=\sevenbf  \def\bf{\fam\bffam\twelvebf}%
\tt  \ttglue=.5em plus.25em minus.15em
\normalbaselineskip=15pt
\setbox\strutbox=\hbox{\vrule height10pt depth5pt width0pt}%
\let\sc=\tenrm  \let\big=\twelvebig  \normalbaselines\rm}
\def\tenpoint{\def\rm{\fam0\tenrm}
\textfont0=\tenrm  \scriptfont0=\sevenrm  \scriptscriptfont0=\fiverm
\textfont1=\teni  \scriptfont1=\seveni  \scriptscriptfont1=\fivei
\textfont2=\tensy  \scriptfont2=\sevensy  \scriptscriptfont2=\fivesy
\textfont3=\tenex  \scriptfont3=\tenex  \scriptscriptfont3=\tenex
\textfont\itfam=\tenit  \def\it{\fam\itfam\tenit}%
\textfont\slfam=\tensl  \def\sl{\fam\slfam\tensl}%
\textfont\ttfam=\tentt  \def\tt{\fam\ttfam\tentt}%
\textfont\bffam=\tenbf  \scriptfont\bffam=\sevenbf
\scriptscriptfont\bffam=\fivebf  \def\bf{\fam\bffam\tenbf}%
\tt  \ttglue=.5em plus.25em minus.15em
\normalbaselineskip=12pt
\setbox\strutbox=\hbox{\vrule height8.5pt depth3.5pt width0pt}%
\let\sc=\eightrm  \let\big=\tenbig  \normalbaselines\rm}
\def\ninepoint{\def\rm{\fam0\ninerm}
\textfont0=\ninerm  \scriptfont0=\sixrm  \scriptscriptfont0=\fiverm
\textfont1=\ninei  \scriptfont1=\sixi  \scriptscriptfont1=\fivei
\textfont2=\ninesy  \scriptfont2=\sixsy  \scriptscriptfont2=\fivesy
\textfont3=\tenex  \scriptfont3=\tenex  \scriptscriptfont3=\tenex
\textfont\itfam=\nineit  \def\it{\fam\itfam\nineit}%
\textfont\slfam=\ninesl  \def\sl{\fam\slfam\ninesl}%
\textfont\ttfam=\ninett  \def\tt{\fam\ttfam\ninett}%
\textfont\bffam=\ninebf  \scriptfont\bffam=\sixbf
\scriptscriptfont\bffam=\fivebf  \def\bf{\fam\bffam\ninebf}%
\tt  \ttglue=.5em plus.25em minus.15em
\normalbaselineskip=11pt
\setbox\strutbox=\hbox{\vrule height8pt depth3pt width0pt}%
\let\sc=\sevenrm  \let\big=\ninebig  \normalbaselines\rm}
\def\eightpoint{\def\rm{\fam0\eightrm}
\textfont0=\eightrm  \scriptfont0=\sixrm  \scriptscriptfont0=\fiverm
\textfont1=\eighti  \scriptfont1=\sixi  \scriptscriptfont1=\fivei
\textfont2=\eightsy  \scriptfont2=\sixsy  \scriptscriptfont2=\fivesy
\textfont3=\tenex  \scriptfont3=\tenex  \scriptscriptfont3=\tenex
\textfont\itfam=\eightit  \def\it{\fam\itfam\eightit}%
\textfont\slfam=\eightsl  \def\sl{\fam\slfam\eightsl}%
\textfont\ttfam=\eighttt  \def\tt{\fam\ttfam\eighttt}%
\textfont\bffam=\eightbf  \scriptfont\bffam=\sixbf
\scriptscriptfont\bffam=\fivebf  \def\bf{\fam\bffam\eightbf}%
\tt  \ttglue=.5em plus.25em minus.15em
\normalbaselineskip=9pt
\setbox\strutbox=\hbox{\vrule height7pt depth2pt width0pt}%
\let\sc=\sixrm  \let\big=\eightbig  \normalbaselines\rm}
\def\twelvebig#1{{\hbox{$\textfont0=\twelverm\textfont2=\twelvesy
        \left#1\vbox to10pt{}\right.\n@space$}}}
\def\tenbig#1{{\hbox{$\left#1\vbox to8.5pt{}\right.\n@space$}}}
\def\ninebig#1{{\hbox{$\textfont0=\tenrm\textfont2=\tensy
        \left#1\vbox to7.25pt{}\right.\n@space$}}}
\def\eightbig#1{{\hbox{$\textfont0=\ninerm\textfont2=\ninesy
        \left#1\vbox to6.5pt{}\right.\n@space$}}}
 
\magnification=\magstephalf\hoffset=0.cm
\voffset=1truecm\hsize=16.5truecm \vsize=21.truecm
\baselineskip=14pt plus0.1pt minus0.1pt \parindent=12pt
\lineskip=4pt\lineskiplimit=0.1pt      \parskip=0.1pt plus1pt
\def\st{\scriptstyle}

\let\st=\scriptstyle


\let\a=\alpha \let\b=\beta   \let\d=\delta  \let\e=\varepsilon
\let\f=\phi  \let\h=\eta      \let\l=\lambda
\let\m=\mu   \let\n=\nu                                  \let\ps=\psi
  \let\s=\sigma \let\t=\tau

\let\D=\Delta \let\F=\Phi  \let\G=\Gamma  \let\L=\Lambda 
\let\O=\Omega


\def\data{
 \number\day/
 \ifcase\month\or January \or February \or March \or
 April \or May \or june \or July \or August \or September
\or October \or November \or December \fi
 /\number\year
 }


\setbox200\hbox{$\scriptscriptstyle \data $}

\newcount\pgn 
\pgn=1
\def\foglio{\veroparagrafo:\number\pgn
\global\advance\pgn by 1}


\global\newcount\numsec
\global\newcount\numfor
\global\newcount\numfig
\global\newcount\numtheo

\gdef\profonditastruttura{\dp\strutbox}

\def\senondefinito#1{\expandafter\ifx\csname#1\endcsname\relax}

\def\SIA #1,#2,#3 {\senondefinito{#1#2}%
   \expandafter\xdef\csname #1#2\endcsname{#3}\else
   \write16{???? ma #1,#2 e' gia' stato definito !!!!}\fi}

\def\etichetta(#1){(\veroparagrafo.\veraformula)
   \SIA e,#1,(\veroparagrafo.\veraformula)
   \global\advance\numfor by 1
   \write15{\string\FU (#1){\equ(#1)}}
   \write16{ EQ \equ(#1) == #1  }}

\def\FU(#1)#2{\SIA fu,#1,#2 }

%
\def\tetichetta(#1){{\veroparagrafo.\verotheo}%
   \SIA theo,#1,{\veroparagrafo.\verotheo}
   \global\advance\numtheo by 1%
   \write15{\string\FUth (#1){\thm[#1]}}%
   \write16{ TH \thm[#1] == #1  }}

\def\FUth(#1)#2{\SIA futh,#1,#2 }
%

\def\getichetta(#1){Fig. \verafigura
 \SIA e,#1,{\verafigura}
 \global\advance\numfig by 1
 \write15{\string\FU (#1){\equ(#1)}}
 \write16{ Fig. \equ(#1) ha simbolo  #1  }}

\newdimen\gwidth

\def\BOZZA{
 \def\alato(##1){
 {\vtop to \profonditastruttura{\baselineskip
 \profonditastruttura\vss
 \rlap{\kern-\hsize\kern-1.2truecm{$\scriptstyle##1$}}}}}
 \def\galato(##1){ \gwidth=\hsize \divide\gwidth by 2
 {\vtop to \profonditastruttura{\baselineskip
 \profonditastruttura\vss
 \rlap{\kern-\gwidth\kern-1.2truecm{$\scriptstyle##1$}}}}}
 \def\talato(##1){\rlap{\sixrm\kern -1.2truecm ##1}}
}

\def\alato(#1){}
\def\galato(#1){}
\def\talato(#1){}

\def\veroparagrafo{\ifnum\numsec<0 A\number-\numsec\else
   \number\numsec\fi}
\def\veraformula{\number\numfor}
\def\verotheo{\number\numtheo}
\def\verafigura{\number\numfig}

\def\Thm[#1]{\tetichetta(#1)}
\def\thf[#1]{\senondefinito{futh#1}$\clubsuit$#1\else
   \csname futh#1\endcsname\fi}
\def\thm[#1]{\senondefinito{theo#1}thf[#1]\else
   \csname theo#1\endcsname\fi}

\def\Eq(#1){\eqno{\etichetta(#1)\alato(#1)}}
\def\eq(#1){\etichetta(#1)\alato(#1)}
\def\eqv(#1){\senondefinito{fu#1}$\clubsuit$#1\else
   \csname fu#1\endcsname\fi}
\def\equ(#1){\senondefinito{e#1}eqv(#1)\else
   \csname e#1\endcsname\fi}


\def\include#1{
\openin13=#1.aux \ifeof13 \relax \else
\input #1.aux \closein13 \fi}
\openin14=\jobname.aux \ifeof14 \relax \else
\input \jobname.aux \closein14 \fi
\openout15=\jobname.aux

\def\fine{\vfill\eject}


%
\newcount\fnts
\fnts=0
\fnts=1 
%
%

\def\smallno{\smallskip\noindent}
\def\medno{\medskip\noindent}
\def\bigno{\bigskip\noindent}

\def\tthsp{\kern .083333 em}

%

\def\indbox#1{\hbox to \parindent{\hfil\ #1\hfil} }

\def\ref[#1]{[#1]}

\def\beginsubsection#1\par{\bigskip\leftline{\it #1}\nobreak\smallskip
            \noindent}

\newfam\msafam
\newfam\msbfam
\newfam\eufmfam
%
%
\ifnum\fnts=0
\def\integer{ { {\rm Z} \mskip -6.6mu {\rm Z} }  }
\def\real{{\rm I\!R}}

\def\Ee{{\rm I\!E}}
\def\Pp{{\rm I\!P}}
\def\mbox{
\vbox{ \hrule width 6pt
   \hbox to 6pt{\vrule\vphantom{k} \hfil\vrule}
   \hrule width 6pt}
}
\def\QED{\ifhmode\unskip\nobreak\fi\quad
  \ifmmode\mbox\else$\mbox$\fi}
\let\restriction=\lceil
%
%
\else
\def\hexnumber#1{%
\ifcase#1 0\or 1\or 2\or 3\or 4\or 5\or 6\or 7\or 8\or
9\or A\or B\or C\or D\or E\or F\fi}
\font\tenmsa=msam10
\font\sevenmsa=msam7
\font\fivemsa=msam5
\textfont\msafam=\tenmsa
\scriptfont\msafam=\sevenmsa
\scriptscriptfont\msafam=\fivemsa
\edef\msafamhexnumber{\hexnumber\msafam}%
\mathchardef\restriction"1\msafamhexnumber16
\mathchardef\square"0\msafamhexnumber03
\def\QED{\ifhmode\unskip\nobreak\fi\quad
  \ifmmode\square\else$\square$\fi}
\font\tenmsb=msbm10
\font\sevenmsb=msbm7
\font\fivemsb=msbm5
\textfont\msbfam=\tenmsb
\scriptfont\msbfam=\sevenmsb
\scriptscriptfont\msbfam=\fivemsb
\def\Bbb#1{\fam\msbfam\relax#1}
\font\teneufm=eufm10
\font\seveneufm=eufm7
\font\fiveeufm=eufm5
\textfont\eufmfam=\teneufm
\scriptfont\eufmfam=\seveneufm
\scriptscriptfont\eufmfam=\fiveeufm

\def\integer{{\Bbb Z}}
\def\natural{{\Bbb N}}
\def\real{{\Bbb R}}

\def\Ee{{\Bbb E}}
\def\Pp{{\Bbb P}}

\fi
%
%
\def\eg{\hbox{\it e.g.\ }}
\def\ie{\hbox{\it i.e.\ }}
\let\sset=\subset
\let\neper=e
\def\nep#1{ \neper^{#1}}
\let\ii=i
\let\emp=\emptyset

\def\ov#1{{1\over#1}}

\def\Pro{\noindent{\it Proof.}\smallno}
\def\Prot#1{\noindent{\it Proof #1.}\smallno}
\def\frac#1#2{{#1 \over #2}}

\def\EE{ \mathop\Ee\nolimits }

\def\sign{\mathop{\rm sign}\nolimits}

\def\diam#1{{\rm diam}\left(#1\right)}
\def\dist#1#2{{\rm dist}\left(#1,#2\right)}

\def\Norm#1{ \left| #1 \right| }
\def\Normm#1{\left\Vert #1 \right\Vert}

\def\normm#1{\Vert #1 \Vert}

\def\scalprod#1#2#3{<#1,#2>_{#3}}
\def\Scalprod#1#2#3{\left<#1,#2\right>_{#3}}
\let\de=\partial

\def\indfn#1{{\bf 1}(#1)}
\def\ptond#1{\left(#1\right)}
\def\pquad#1{\left[#1\right]}
\def\pgraf#1{\left\{#1\right\}}

\def\qed{\QED\medno}

\def\cE{{\cal E}}

\def\cG{{\cal G}}
\def\cI{{\cal I}}
\def\cS{{\cal S}}

\outer\def\nproclaim#1 [#2]#3. #4\par{\medbreak \noindent
   \talato(#2){\bf #1 \Thm[#2]#3.\enspace }%
   {\sl #4\par }\ifdim \lastskip <\medskipamount 
   \removelastskip \penalty 55\medskip \fi}

\def\thmm[#1]{#1}
\def\teo[#1]{#1}
%
%
\def\sttilde#1{%
\dimen2=\fontdimen5\textfont0
\setbox0=\hbox{$\mathchar"7E$}
\setbox1=\hbox{$\scriptstyle #1$}
\dimen0=\wd0
\dimen1=\wd1
\advance\dimen1 by -\dimen0
\divide\dimen1 by 2
\vbox{\offinterlineskip%
   \moveright\dimen1 \box0 \kern - \dimen2\box1}
}
%

\def\inda{_L^{L/2}}
\def\indb{_L^M}
\def\indc{_L^{\infty}}

\def\dual{\left(\integer^2\right)^*}
\def\FF{{\bf \F}}
\def\cmax{c_{\rm max}}
\def\omean{\mathop{{\bf\overline E}_L^M}\nolimits}
\def\oce{\mathop{\overline \cE_L^M}\nolimits}
\def\bnu{\bar\nu_L^M}
\def\bnl{\bar\nu_L}
\def\oml{\mathop{{\bf\overline E}_L}\nolimits}
\def\dom#1{{\rm Dom}(#1)}
%
%
\font\ttlfnt=cmcsc10 scaled 1200 
%
\begingroup
\nopagenumbers
\footline={}
%
%
\def\author#1
{\vskip 18pt\tolerance=10000
\noindent\centerline{\caps #1}\vskip 1truecm}
%
%
%
%
\def\abstract#1
{
\noindent{\bf Abstract.\ }#1\par}


\vskip 1cm
\centerline{\ttlfnt Equilibrium Fluctuations}
\centerline{\ttlfnt for a One-Dimensional Interface}
\centerline{\ttlfnt in the Solid on Solid Approximation}
\vskip 0.5truecm
\author{Gustavo Posta$^*$}
\abstract{\ninerm
An unbounded one-dimensional solid-on-solid model with integer heights is studied.
Unbounded here means that there is no {\it a priori} restrictions on the discrete gradient of the interface.
The interaction Hamiltonian of the interface is given by a finite range part, proportional to the sum of height differences, plus a part of exponentially decaying long range potentials.
The evolution of the interface is a reversible Markov process.
We prove that if this system is started in the center of a box of size $L$ after a time of order $L^3$ it reaches, with a very large probability, the top or the bottom of the box.
}

\vfill



{\parindent=0pt
\footnote{}{$^*$ Dipartimento di Matematica F.~Brioschi, 
   Politecnico di Milano, Piazza Leonardo da Vinci 32, I-20133 Milano, Italy. e-mail: {\tt gustavo.posta@polimi.it}}

}

\fine

\endgroup


\expandafter
 \ifx\csname sezioniseparate\endcsname\relax
  \fi

\numsec=1
\numfor=1
\numtheo=1
\pgn=1

\beginsection 1. Introduction

The rigorous analysis of Glauber dynamics for classical spin systems 
when the inverse temperature $\b$ is such that the static system does not 
undergo phase transition in the thermodynamic limit has been in the last 
years the argument of many important works.
In particular we refer to [12], [8], [9] and references in these papers.

A natural question is what happens when the thermodynamic parameters
are such that there is a phase transition.  To be concrete consider
the stochastic ferromagnetic Ising model in $\L_L\equiv [1,L] \times
[1,L] \cap \integer ^2$ in absence of external field with free
boundary conditions and let us suppose that the inverse temperature
$\b$ is much larger than the critical one.  Then, as well known (see
[4]), any associated infinite volume (\ie $L=+\infty$) Glauber
dynamics is not ergodic.  However, for $L<+\infty$, every associated
Glauber dynamics is ergodic, because it is an irreducible finite-state
Markov chain.  The problem of how this system behaves at the
equilibrium is discussed in [6].  In that paper is proved that the
system fluctuates between the ``$+$'' phase and the ``$-$'' on two
distinct time scales.  In a first time interval the system creates a
layer of the opposite phase, separated from the initial one by a
one-dimensional interface. In a second time interval the interface
moves until the new phase invades the system.  The time the system
spends for this phase transition is exponential in $L$.  In this paper
we will study the motion of the interface in the solid-on-solid (SOS)
approximation.  In particular we will show that after the formation of
the interface, a time of order $L^3$ is sufficient to reach the
opposite phase.

The SOS model studied here is a one-dimensional random interface (or surface) with integer heights.
There is no restriction on the discrete gradient of the surface (unbounded SOS),
thus the configuration space is $\integer^L$.
The interaction Hamiltonian of the interface is given by the usual energy proportional to the sum of height differences plus a part of exponentially decaying potential which mimics the long-range dependence that the interface feels due to the surrounding bulk phases.
This long-range interaction is small in the regime studied, but the potentials involved are of unbounded range, so it is not completely trivial to handle.
We restrict the interface to stay in a finite box $[1,L]\times[-M,M]$.
This gives (see Section~2 for more details) the Gibbs measure on $\integer^L$
$$
 \m\indb(\f)\equiv
 \frac{\indfn{\normm{\f}_{\infty}\leq M}\nep{-\b H_L(\f)-W_L^M(\b,\f)}}
 {Z\indb(\b)}.
$$
The evolution of the interface is described by a reversible Markov process with generator
$$
 (G\indb f)(\f)\equiv\sum_\ps c\indb(\f,\ps)\pquad{f(\ps)-f(\f)},
$$
where the jump rates are bounded and such that only transitions of the form $(\phi_1,\ldots,\phi_L)\mapsto(\phi_1,\ldots,\phi_k\pm 1,\ldots\phi_L)$, for some $k=1,\ldots,L$, are allowed.
We study this process for $M=L/2$, \ie in a ``box'' of size $L$ and we prove, in a sense given precisely by Theorem~\thf[v0], that if the interface is started in the center of the box then in a time of order $L^3$ it reaches the bottom or the top of the box. 

The technique we use to get this result is a mix of analytical and probabilistic tools.
In fact it is standard to obtain estimates on the exit time distribution of a reversible Markov process from a region if one has an estimate from below of the spectral gap of its generator with Dirichlet boundary conditions.
We will not give bounds on the spectral gap of $G_L^{L/2}$ but on the generator of an auxiliary simpler process and we conclude with a coupling argument.

This paper completes the study of the asymptotic properties of the solid-on-solid model started
in [10] where a similar model, in the presence of boundary conditions, is studied.

\bigno
{\bf Acknowledgments:} I would like to thank Fabio Martinelli who posed this problem to me and helped me with many constructive discussions.
\fine


\expandafter
 \ifx\csname sezioniseparate\endcsname\relax
  \fi

\numsec=2
\numfor=1
\numtheo=1
\pgn=1

\beginsection 2. Notation and Results

Our sample space is
$\O_L\equiv\integer^L$ for fixed $L\in\natural$.
Configurations, \ie
elements of the sample space will be denoted by Greek letters,
\eg $\f=(\f_1,\ldots,\f_L)\in\O_L$.
Given $M\in\natural\cup\pgraf{+\infty}$, and $\b>0$
one defines the {\it energy\/} associated with the configuration $\f\in\O_L$ as:
$$
 \b H_L(\f)+W_L^M(\b,\f).
$$
Here $H_L$ is a local interaction:
$$
 H_L(\f)\equiv\sum_{k=1}^{L-1}|\f_{k+1}-\f_k|,
\Eq(1.s.4.0.1)
$$
while $W_L^M(\b,\f)$ is a long range interaction that will be defined below.

Consider the lattices $\integer^2$ and $\dual\equiv(1/2,1/2)+\integer^2$ 
as graphs embedded in $\real^2$ equipped with the usual Euclidean metric
denoted with $\dist{\cdot}{\cdot}$:
$$
        \dist{(x_1,y_1)}{(x_2,y_2)}\equiv \sqrt{(x_1-x_2)^2+(y_1-y_2)^2}.
$$
For every $\f\in\O_L$ the {\it contour associated with} $\f$ is the 
subset of $\real^2$ defined by:
$$
 \G(\f)
  \equiv\Big[\bigcup_{i=1}^L\big\{(x,y):x\in(i-1,i),\,y=\f_i\big\}\Big]
  \cup\Big[\bigcup_{i=1}^{L-1}\big\{(x,y):x=i,\,y\in[\f_i\land\f_{i+1},\f_i\lor\f_{i+1}]\big\}\Big].
$$
For $A\subset\real^2$ and $p\in\real^2$, the distance of $x$ from $A$ is defined as $\dist{p}{A}\equiv\inf\{\dist{x}{y}:y\in A\}$.
If we denote with
$$
 \D(\f)\equiv\pgraf{p^*\in\ptond{\integer^2}^*:
  \dist{p^*}{\G(\f)}=\ov2,\hbox{ or }\dist{p^*}{\G(\f)}=\ov{\sqrt{2}}}
$$
the set of sites {\it attached\/} to the contour $\G(\f)$ and define
$$
 V_L^M\equiv\ptond{[-1/2,L+1/2]\times[-(M+1)/2,(M+1)/2]}\cap\dual,
\Eq(2.defV)
$$
the long range interaction may be written as:
$$
 W_L^M(\b,\f)\equiv\sum_{\st\L\cap\D(\f)\not=\emptyset\atop\st\L\sset V_L^M}\F(\b,\L).
$$
The sum is over all $\L\sset V_L^M$ connected in the sense of the dual graph $\dual$.
The {\it potential} $\F(\cdot,\cdot)$ is a function satisfying (see. [2]):
\item{$i)$} there exists $\bar\b>0$ such that for every $\b>\bar\b$ we have:
 $$
  \sum_{\st \L\ni p^*\atop\st \diam{\L}\geq k}|\F(\b,\L)|\leq\nep{-m(\b)k}
 \Eq(100)
 $$
for any $k>0$ and $p^*\in\dual$.
Here:
\itemitem{$a)$} the sum is over all $\L\sset\dual$ connected and such 
that $\L\ni p^*$;
\itemitem{$b)$} $m(\b)$ is a positive function such that $m(\b)\to+\infty$ for $\b\to+\infty$;
\itemitem{$c)$} $\diam{\L}$  is the diameter of the set $\L$.
\item{\it ii)} For every $p^*\in\dual$
$$
 \F(\b,\L+p^*)=\F(\b,\L).
  \Eq(dec2)
$$

\noindent
It is easy to check that for every $\b>0$ and $M\in(0,+\infty]$ there 
exists finite
the {\it partition function}:
$$
 Z\indb(\b)\equiv\sum_{\st\f\in\O_L}\indfn{\normm{\f}_{\infty}\leq M}
   \nep{-\b H_L(\f)-W_L^M(\b,\f)},
\Eq(1.s.4.0.2)
$$
where $\normm{\f}_{\infty}\equiv\max_{1\leq i\leq L}|\f_i|$.
Thus one can define on $\O_L$  the {\it Gibbs measure\/}:
$$
 \m\indb(\f)\equiv
 \frac{\indfn{\normm{\f}_{\infty}\leq M}\nep{-\b H_L(\f)-W_L^M(\b,\f)}}
 {Z\indb(\b)}.
\Eq(1.s.2.0)
$$
This measure represents the equilibrium of the system.
The dynamics of the system is a continuous time Markov chain with values in $\O_L$
and stationary measure $\m\indb$.
This process will be defined by means of its generator.

For every $\f\in\O_L$, $k=1,\ldots L$, define:
$$
 \f\pm\d_k=(\f_1,\ldots,\f_k\pm1,\ldots,\f_L),
$$
and the {\it jump rates}:
$$
 c\indb(\f,\ps)=
 \cases{
  \ptond{\m\indb(\ps)\over\m\indb(\f)}^{\ov2} 
  &if $\m\indb(\f)>0$ and $\ps=\f\pm\d_k$ for some $k=1,\ldots,L$;\cr
   & \cr
  0&otherwise.\cr
 }
\Eq(defrates2)
$$
It is simple to prove that there exists a unique Markov process $\FF\equiv\pgraf{\FF(t):t\geq0}$
with generator
$$
 (G\indb f)(\f)\equiv\sum_\ps c\indb(\f,\ps)\pquad{f(\ps)-f(\f)}.
\Eq(defg)
$$
Moreover $G\indb$ is self-adjoint in $L^2(\m\indb)$, \ie $\FF$ is {\it reversible} and
$G\indb$ is negative semidefinite.
The absolute value of the largest negative eigenvalue of $G\indb$ is denoted by $\l_1(G\indb)$ and it is called
{\it spectral gap\/} of $G\indb$.

We will give a direct construction of $\FF$ in Section~4.
More precisely (see Proposition~\thf[v6]) we will define a measurable space
$({\bf\O}_L,{\bf F}_L)$ and a family of probability measures on it 
$\pgraf{\Pp_\f:\f\in\O_L}$ such that:
\item{$i)$}
for every $\f\in\O_L$
the process $\FF$ is a Markov process on $({\bf\O}_L,{\bf F}_L,\Pp_\f)$
with generator $G_L^{\emp,M,W}$;
\item{$ii)$}
 $\Pp_\f(\FF(0)=\f)=1$.

\noindent
For every measurable set $A\sset\O_L$ define the {\it first exit time\/} from $A$ as:
$$
 \t(A^c)\equiv\inf\pgraf{t\geq 0:\FF(t)\in A^c}.
$$
The main result of this paper is:

\nproclaim Theorem [v0].
Let $\FF$ be the process associated with the generator $G\inda$. Fix $\a\in\ptond{0,1/4}$, $\e\in\ptond{0,1/100}$
and define
$A\equiv\pgraf{\f\in\O_L:\normm{\f}_{\infty}\leq (1-\e)L/2}$.
Then there exists $\bar\b>0$ and for every $\b>\bar\b$ constants $K_1(\b)$, $K_2(\b)$, $K_3(\b)$ and $K_4(\b)>0$
such that:
$$
 \int_{\O_L}d\m_L^{\emp,L/2,W}(\f|\,\normm{\f}_{\infty}\leq\a L)\Pp_\f(\t(A^c)>t)
  \leq \a^{-1} K_1\nep{-K_2\left({t\over L^3}\land L\right)}+\a^{-1} K_3\nep{-K_4\e L}
\Eq(v0.1)
$$
for any $L>0$.

This result can be read in the following way: starting the interface in a square box of size $L$, from an initial condition randomly chosen under $\m_L^{\emp,L/2,W}(\cdot|\,\normm{\cdot}_{\infty}\leq\a L)$ (\ie the interface is forced to stay at least $\a L$ far away from the top or the bottom of the box), the probability of reaching within $\e L$ of the top or the bottom of the box in time bigger than $t$ is exponentially small in $t/L^3$.

\medskip\noindent
{\bf Remark 2.2.}
\ In what follows we use constants $K_1,K_2,\dots$ in the statement of theorems, propositions and so on, while we use constants $C_1,C_2,\ldots$ inside proofs.
The reader should be warned that the use of constants is coherent only inside the same structure. This means, \eg, that constants which appears in the proof or in the statement of a proposition may differ from constants, with the same name, which appears in the proof or in the statement of a different proposition.   

\fine


\expandafter
 \ifx\csname sezioniseparate\endcsname\relax
  \fi

\numsec=3
\numfor=1
\numtheo=1
\pgn=1

\beginsection 3. Proof of Main Result

In this section we will prove our main result Theorem~\thm[v0].
The technique we will use is the following.
We can estimate the first exit time of a reversible Markov
process from a region $A$ by bounding from below the spectral gap of the generator of the process.
Actually we will not bound the spectral gap of the process $\FF$.
Instead we will estimate the spectral gap of a simpler auxiliary process $\bar\FF$ 
that in the region $A$ is similar to $\FF$.
Then a {\it coupling\/} argument (Proposition~\thf[v6]) concludes the proof.
\bigskip
We begin this section introducing the auxiliary process above mentioned.
Fix $L\in\natural$, $\b$ and $M>0$ and
define on $\O_L$ the probability measure
$$
 \bar\m\indb(\f)\equiv\indfn{|\f_1|\leq M}
  {\nep{-\b H_L(\f)- W\indc(\b,\f)}\over \bar Z\indb(\b)},
\Eq(v0.2)
$$
where
$$
 \bar Z\indb(\b)\equiv\sum_\f\indfn{|\f_1|\leq M}
  \nep{-\b H_L(\f)- W\indc(\b,\f)}
$$
and (see Section~2)
$$
 W\indc(\b,\f)\equiv\sum_{\L\cap\D(\f)\not=\emp\atop\L\sset V\indc}\F(\b,\L).
$$
The process $\bar\FF$ is defined by means of its generator on $L^2(\O_L,\bar\m\indb)$
$$
 (\bar G\indb f)(\f)\equiv\sum_\ps \bar c\indb(\f,\ps)\pquad{f(\ps)-f(\f)},
\Eq(defgW)
$$
where
$$
 \bar c\indb(\f,\ps)\equiv
 \cases{
  \ptond{\bar\m\indb(\ps)\over\bar\m\indb(\f)}^{\ov2} 
   &if $\bar\m\indb(\f)>0$ and $\ps=\f\pm\d_k$ for some $k=1,\ldots,L$;\cr
   & \cr
  0&otherwise.\cr
 }
\Eq(defrates3)
$$
It simple to check that $\bar G\indb$ is a self-adjoint Markov generator 
which defines a unique Markov process.
Because for $\f\in A$, defined in Theorem~\thm[v0], the jump rates of $\bar G\inda$ are very close to the jump rates 
of $G\inda$  (see \equ(defrates2)), the processes $\bar\FF$ and $\FF$ evolve in a similar way as long as they 
remains within $A$.
This fact is formally proved in the following result which gives also a 
direct construction of the processes.

\nproclaim Proposition [v6].
It is possible to construct a family of probability spaces $({\bf\O}_L,{\bf F}_L,\Pp_{\f,\bar\f})$ and a process $\{(\FF(t),\bar\FF(t)):t\geq 0\}$ taking values in $\O_L\times\O_L$,  with  $\Pp_{\f,\bar\f}(\FF(0)=\f,\ \bar\FF(0)=\bar\f)=1$ for any $(\f,\bar\f)\in\O_L\times\O_L$, and such that:
\item{$i)$} $\FF$ and $\bar\FF$ are Markov processes with generators $G\inda$ and $\bar G\inda$ respectively;
 \item{$ii)$} if we define $\s\equiv\inf\pgraf{t\geq0:\FF(t)\not=\bar\FF(t)}$ and  $\bar\t(A^c)\equiv\inf\pgraf{t\geq0:\bar\FF(t)\in A^c}$, where $A$ is defined in Theorem~\thm[v0], then there exist $K_1(\b)$ and $K_2(\b)>0$, with $K_2(\b)\to+\infty$ for $\b\to+\infty$,
 such that for every $\f\in A$:
 $$
  \Pp_{\f,\f}(\s\leq t,\ \s\leq\bar\t)\leq K_1 Lt \nep{-K_2\e L}.
 \Eq(v6.1)
 $$

\noindent
This proposition will be proved in Section~4.
 
The advantage in considering the process $\bar \FF$ instead of $\FF$ is 
that the first one is simpler to study because it has no interaction 
with the top and the bottom of the box $V\indb$.
In particular in Section~5 we will prove the following result on the 
spectral gap $\l_1(\bar G\indb)$:

\nproclaim Proposition [v5].
There exists $\bar\b>0$ such that for every $\b\geq\bar\b$ there exists $K_1(\b)>0$,
so that:
$$
 \l_1(\bar G\indb)\geq K_1L^{-1}(L^{-1}\land M^{-2})
\Eq(v5.1)
$$
for every $L$ and $M>0$.

\noindent
We are now in a position to prove Theorem~\thm[v0].

\Prot{of Theorem~\thm[v0]}
In this proof we simplify notation writing $G\equiv G\inda$,
$\bar G\equiv \bar G\inda$, $\m\equiv \m\inda$,
$\bar\m\equiv \bar\m\inda$ and
$$
 \eqalign{
  \t\equiv\t(A^c)&=\inf\pgraf{t\geq0:\normm{\FF(t)}_{\infty}>(1-\e)L/2},\cr
  \bar\t\equiv\bar\t(A^c)&=\inf\pgraf{t\geq0:\normm{\bar\FF(t)}_{\infty}>(1-\e)L/2}.\cr
 }
$$
Recall that $\s$ was defined in Proposition~\thm[v6] as the first time 
such that $\bar\FF(\s)\not=\FF(\s)$ and
suppose that $\f\in A$. Then for any $t>0$
$$
 \Pp_{\f,\f}(\t>t)
 =\Pp_{\f,\f}(\t>t,\ \t=\bar\t)+\Pp_{\f,\f}(\t>t,\ \t\not=\bar\t)
 \leq\Pp_{\f,\f}(\bar\t>t)+\Pp_{\f,\f}(\s\leq\bar\t)
$$
and any $s>0$
$$
 \eqalign{
  &\quad\Pp_{\f,\f}(\s\leq\bar\t)
   =\Pp_{\f,\f}(\s\leq\bar\t,\ \s\leq s )+\Pp_{\f,\f}(\s\leq\bar\t,\ \s>s)\leq\cr
  &\leq\Pp_{\f,\f}(\s\leq\bar\t,\ \s\leq s )+\Pp_{\f,\f}(\bar\t> s).\cr
 }
$$
In conclusion for any $t,s>0$
$$
 \Pp_{\f,\f}(\t>t)
 \leq2\Pp_{\f,\f}(\bar\t>t\land s)+\Pp_{\f,\f}(\s\leq\bar\t,\ \s\leq s).
\Eq(v0.3)
$$
Define $B\equiv\pgraf{\f\in\O_L:\normm{\f}_{\infty}\leq\a L}$,
where $\a\in(0,1/4)$.
By \equ(v0.3) we obtain:
$$
 \int_{\O_L}d\m(\f|B)\Pp_{\f,\f}(\t>t)
 \leq 2\int_{\O_L}d\m(\f|B)\Pp_{\f,\f}(\bar\t>t\land s)
   +\sup_{\f\in B}\Pp_{\f,\f}(\s\leq s).
\Eq(v0.4)
$$
We are going to bound from above
 the first term on the right hand side of \equ(v0.4).
This is done using a Markov process with killing (see [11]).
The {\it Dirichlet form\/} associated with the generator $\bar G$ is the 
positive-semidefinite bilinear form
$$
 \bar\cG(f,g)\equiv-\scalprod{\bar G f}{g}{L^2(\O_L,\bar\m)},
$$
where $\scalprod{\cdot}{\cdot}{L^2(\O_L,\bar\m)}$ is the scalar product in $L^2(\O_L,\bar\m)$.   
Because $\bar\m(A)>0$ one can define the positive 
semidefinite bilinear form
$$
 \bar\cG_A(f,g)\equiv\bar\cG(\hat f,\hat g)
\Eq(2.1)
$$
with form domain
$$
 \dom{\bar\cG_A}\equiv\pgraf{f\in L^2(A,\bar\m):\hat f\in\dom{\bar\cG}}
$$
where
$$
 \hat f(\f)=\cases{f(\f) &if $\f\in A$;\cr 0 & otherwise.\cr}
$$
Standard functional analysis methods shows that there exists a unique 
positive-semidefinite self-adjoint (in $L^2(A,\bar\m)$) operator $\bar G_A$ such that:
$$
 \bar\cG_A(f,g)\equiv-\scalprod{\bar G_A f}{g}{L^2(A,\m)}.
$$
This operator has a probabilistic interpretation, it is the generator of 
a process which evolves according to $\bar G$ as long as it stays within $A$, 
but is killed when it tries to jump outside $A$ (see [11]).
The semigroup $\nep{t\bar G_A}$ generated by $\bar G_A$ is {\it sub-stochastic\/} and 
we have:
$$
 (\nep{t\bar G_A}f)(\f)=\EE_{\f,\f}\pquad{\indfn{\bar\t>t}f(\bar\FF(t))}
$$
for every $f\in L^2(A,\bar\m)$.
In particular taking $f\equiv{\bf1}$ we obtain:
$$
 \Pp_{\f,\f}(\bar\t>t)=(\nep{t\bar G_A}{\bf 1})(\f).
$$
Thus:
$$
 \eqalign{
  &\quad\int_{\O_L}d\m(\f|B)\Pp_{\f,\f}(\bar\t>t)
   =\Scalprod{(\nep{t\bar G_A}{\bf 1})}{{d\m(\cdot|B)\over d\bar\m}}{L^2(A,\bar\m)}\leq\cr
  &\leq\Normm{(\nep{t\bar G_A}{\bf 1})}_{L^2(A,\bar\m)}
    \Normm{{d\m(\cdot|B)\over d\bar\m}}_{L^2(A,\bar\m)}
   \leq\Normm{{d\m(\cdot|B)\over d\bar\m}}_{L^2(A,\bar\m)}\nep{-\l_1(\bar G_A)t}.\cr
 }
\Eq(v0.5)
$$
Spectral theorem has been used in the last line.
The spectral gap $\l_1(\bar G_A)$ is characterized by the following 
variational property:
$$
 \l_1(\bar G_A)
 =\inf_{f\in L^2(A,\bar\m)\atop f\perp{\bf1}}
  \frac{\bar\cG_A(f,f)}{\normm{f}_{L^2(A,\bar\m)}^2}.
$$
This relation and \equ(2.1) imply:
$$
 \l_1(\bar G_A)
 =\inf_{f\in L^2(A,\bar\m)\atop f\perp{\bf1}}
  \frac{\bar\cG(\hat f,\hat f)}{\normm{\hat f}_{L^2(\O_L,\bar\m)}^2}
 \geq\inf_{f\in L^2(\O_L,\bar\m)\atop f\perp{\bf1}}
  \frac{\bar\cG(f,f)}{\normm{f}_{L^2(\O_L,\bar\m)}^2}
 =\l_1(\bar G).
$$
By Proposition~\thm[v5] we know that $\l_1(\bar G)\geq C_1(\b) L^{-3}$, so 
\equ(v0.5) yields:
$$
 \int_{\O_L}d\m(\f|B)\Pp_{\f,\f}(\bar\t>t)
  \leq\nep{-{tC_1\over L^3}}\Normm{{d\m(\cdot|B)\over d\bar\m}}_{L^2(A,\m)}.
$$
We claim that there exists $C_2(\b)>0$ such that: 
$$
 \Normm{{d\m(\cdot|B)\over d\bar\m}}_{L^2(A,\m)}\leq\a^{-1}C_2,
$$
this simple technical bound is proved in the appendix (see Lemma~\thf[v8]).
In conclusion:
$$
 \int_{\O_L}d\m(\f|B)\Pp_{\f,\f}(\bar\t>t)
  \leq\a^{-1}C_2\nep{-{C_1t\over L^3}}
$$
and by \equ(v0.4)
$$
 \int_{\O_L}d\m(\f|B)\Pp_{\f,\f}(\t>t)
  \leq\a^{-1}C_2\nep{-{C_1(t\land s)\over L^3}}+
   \sup_{\f\in B}\Pp_{\f,\f}(\s\leq s).
$$
Taking $s=L^4$  by \equ(v6.1) we obtain:
$$
 \int_{\O_L}d\m(\f|B)\Pp_{\f,\f}(\t>t)
  \leq C_2\a^{-1}\nep{-C_1\left({t \over L^3}\land L\right)}+C_3\a^{-1}\nep{-C_4\e L},
$$
\ie  \equ(v0.1).
\qed
\fine


\expandafter
 \ifx\csname sezioniseparate\endcsname\relax
  \fi

\numsec=4
\numfor=1
\numtheo=1
\pgn=1

\beginsection 4. The Coupling

In this section we will construct explicitly a stochastic coupling between $\FF$ and $\bar\FF$,
in particular we will prove Proposition~\thm[v6].
The technique we use is an application of the so called {\it basic coupling}.
This is a coupling between jump processes such that the processes jump together
as long as possible, considering the constraint they have to jump 
with their own jump rates.
Because the jump rates of $\FF$ and $\bar\FF$ are very close, when they are in $A=\pgraf{\f\in\O_L:\normm{\f}_{\infty}\leq (1-\e)L/2}$, the two 
processes will evolve identically for a long time.
\bigskip
The first step in the construction of the coupling is to show that the 
jump rates of  
$\FF$ and $\bar\FF$ are close in $A$.

\nproclaim Lemma [v7].
For any $\b>\bar \b$ there exists $K_1(\b)$ and $K_2(\b)>0$ such that:
$$
 \sup_{\st k=1,\ldots,L\atop\st\f\in A}
  \Norm{{\bar c\inda(\f,\f\pm\d_k)\over c\inda(\f,\f\pm\d_k)}-1}
  \leq K_1\nep{-m(\b)\e K_2 L}
\Eq(v7.1)
$$
for every $L>0$.

\Pro
To simplify we adopt the notation of the last section an we write
$\bar c\equiv\bar c\inda$, $c\equiv c\inda$ and $\bar\m\equiv \bar\m\inda$.
Notice that
$$
 \pquad{{\bar c(\f,\f+\d_k)\over c(\f,\f+\d_k)}}^2=
 {\bar\m(\f+\d_k)\m(\f)\over\bar\m(\f)\m(\f+\d_k)}
  =\nep{[(\de_k^+W\inda)(\b,\f)-(\de_k^+W\indc)(\b,\f)]}
\Eq(v7.2)
$$
for every $\f\in A$.
Here and later $(\de_k^\pm f)(\f)\equiv f(\f\pm\d_k)-f(\f)$.
If we define $p_k(\f)\equiv(k,\f_k)$ and $I_k\equiv\pgraf{p\in\dual:\dist{p}{p_k}\leq 4}$
it simple to prove that:
$$
 |(\de_k^+W\inda)(\b,\f)-(\de_k^+W\indc)(\b,\f)|
  \leq\sum_{\st\L\cap I_k\not=\emp\atop\st\L\cap (V\inda)^c\not=\emp}|\F(\b,\L)|.
\Eq(v7.3)
$$
If $\f\in A$, $\L\cap I_k\not=\emp$ and $\L\cap (V\inda)^c\not=\emp$ then there exists a constant $C_1>0$ such that
$\diam{\L}\geq\e C_1 L$.
We can use the condition \equ(100) to bound the sum on the right hand 
side of \equ(v7.3).
This gives
$$
 |(\de_k^+ W_L^{L/2})(\f)-(\de_k^+W_L^\infty)(\f)|\leq 16\nep{-m(\b)\e C_1L}.
$$
From this relation and \equ(v7.2) we have \equ(v7.1).
\qed

We can now prove the main result of this section.

\Prot{of Proposition \thm[v6]}
We use the basic coupling.
To any site $k=1,\ldots,L$ we associate two independent Poisson 
processes, each one with rate
$\cmax\equiv\sup_{\f,\ps}\pgraf{c(\f,\ps),\bar c(\f,\ps)}$.
We will denote these processes $\{N_{k,t}^+:t\geq0\}$ and $\{N_{k,t}^-:t\geq0\}$
while the arrival times of each process are denoted by $\{\t_{k,n}^+:n\in\natural\}$ and
$\{\t_{k,n}^-:n\in\natural\}$ respectively.
Assume that the Poisson processes associated to different sites are also 
mutually independent.
We say that at each point in the space-time of the form
$(k,\t_{k,n}^+)$ there is a ``$+$'' mark and that at each point of the form  $(k,\t_{k,n}^-)$
there is a ``$-$'' mark.

Next we associate to each arrival time $\t_{k,n}^\pm$ a random variable
$U_{k,n}^\pm$ with uniform distribution in $[0,1]$.
All these random variables are assumed to be independent among 
themselves and independent from the previously introduced Poisson processes.
Obviously there exists a probability space such that all these objects are defined.
We have to say now how the various processes are constructed on this space.
The process
$\FF$ (resp. $\bar\FF$) is defined in the following manner.
We know that almost surely the arrival times $\t_{k,n}^*$, $k=1,\ldots,L,\ n\in\natural,\ *=\pm$
are all distinct.
We update the state of the process each time there is a mark at some $k=1,\ldots,L$
according to the following rule.
\item{$\bullet$}
 If the mark that we are considering is at the point $(k,\t_{k,n}^*)$, 
 with $*=\pm$,
 and the configuration of $\FF$ (resp. of $\bar\FF$) immediately before time $\t_{k,n}^*$ was $\f$,
 (resp. $\bar\f$) then the configuration immediately after $\t_{k,n}^*$ of $\FF$ (resp. of
 $\bar\FF$) will be $\f\pm\d_k$ (resp. $\bar\f\pm\d_k$) if an only if
 $$
  c(\f,\f\pm\d_k)>U_{k,n}^\pm\cmax
   \qquad\qquad\hbox{(resp. $\bar c(\f,\f\pm\d_k)>U_{k,n}^\pm\cmax$)}.
 $$
 Else the configuration remains the same.

\noindent
It is easy to check that this construction satisfies condition $i)$  of the proposition.
It remains to prove \equ(v6.1).

Define $N_t=\sum_{k=1}^L(N_{k,t}^++N_{k,t}^-)$.
This process counts the number of possible updating
of the processes $\FF$ and $\bar\FF$ in the interval $[0,t]$.
It is clear that $\pgraf{N_t:t\geq0}$ is a Poisson process with rate 
$\l\equiv2L\cmax$.
For $\f\in A$ we have:
$$
 \Pp_{\f,\f}(\s\leq t,\ \s\leq\bar\t)
 =\nep{-\l t}\sum_{n=1}^{+\infty}{(\l t)^n\over n!}\Pp_{\f,\f}(\s\leq t,\ \s\leq\bar\t|N_t=n).
\Eq(v6.2)
$$
To bound from above $\Pp_{\f,\f}(\s\leq t,\ \s\leq\bar\t|N_t=n)$
observe that if $\FF$ and $\bar\FF$ are initially in the same state
$\f\in A$ and if  $N_t=n$, \ie there were $n$ possible updating in $[0,t]$,
then it possible that $\s\in[0,t]$ if and only if for some $i=1,\ldots,n$, $\ps\in\O_L$ and $k=1,\ldots,L$ happens that:
$$
 \bar c(\ps,\ps\pm\d_k)> U_{k,i}^\pm\cmax
 \qquad\hbox{and}\qquad
 c(\ps,\ps\pm\d_k)\leq U_{k,i}^\pm\cmax
$$
or
$$
 \bar c(\ps,\ps\pm\d_k)\leq U_{k,i}^\pm\cmax
        \qquad\hbox{and}\qquad
        c(\ps,\ps\pm\d_k)> U_{k,i}^\pm\cmax.
$$
The probability of this event, for fixed $i=1,\ldots,n$, $\ps\in\O_L$ and $k=1,\ldots,L$, is:
$$
 \Norm{{\bar c(\ps,\ps\pm\d_k)\over\cmax}-{c(\ps,\ps\pm\d_k)\over\cmax}}.
$$
Moreover because $\s\leq\bar\t$ we have:
$$
 \Pp_{\f,\f}(\s\leq t,\ \s\leq\bar\t|N_t=n)
  \leq1-\bigg(1-\sup_{\st k=1,\ldots,L\atop\st\ps\in A}
   \Norm{{\bar c(\ps,\ps\pm\d_k)\over\cmax}-{c(\ps,\ps\pm\d_k)\over\cmax}}\bigg)^n.
$$
By Lemma~\thf[v7] we have
$$
 \sup_{\st k=1,\ldots,L\atop\st\ps\in A}
  \Norm{{\bar c(\ps,\ps\pm\d_k)\over\cmax}-{c(\ps,\ps\pm\d_k)\over\cmax}}\leq C_1\nep{-m(\b)C_2\e L}
$$
which gives
$$
 \Pp_{\f,\f}(\s\leq t,\ \s\leq\bar\t|N_t=n)\leq1-\pquad{1-C_1\nep{-m(\b)C_2\e L}}^n.
$$
This estimate together with \equ(v6.2) gives \equ(v6.1).
\qed
\fine


\expandafter
 \ifx\csname sezioniseparate\endcsname\relax
  \fi

\numsec=5
\numfor=1
\numtheo=1
\pgn=1

\beginsection 5. Spectral Gap for $\bar G\indb$

In this section we will prove Proposition~\thf[v5].
The strategy of the proof is the following.
By a simple change of variables the {\it Glauber\/} dynamics associated 
with $\bar G\indb$ becomes a {\it Kawasaki\/} type dynamics, while the 
measure $\bar\m\indb$ becomes, in the new variables, a product measure 
perturbed by an infinite range interaction.
This interaction term is small if $\b$ is large.
Without this perturbation term the result is very simple to prove. 
The presence of this extra term requires a little extra work.
\bigskip
Consider the random variables defined by
$$
 \eqalign{
  \h_1(\f)&\equiv\f_1\cr
  \h_2(\f)&\equiv\f_2-\f_1\cr
  &\ldots\cr
  \h_{L}(\f)&\equiv\f_{L}-\f_{L-1}.\cr
 }
\Eq(v.defeta)
$$
Obviously the map $\O_L\ni\f\mapsto\O_L\in\h$ is bijective with inverse 
map
$$
  T:\O_L\ni(\h_1,\ldots,\h_L)\mapsto(\h_1,\h_1+\h_2,\ldots,\h_1+\ldots+\h_{L})\in\O_{L}.
$$
The distribution of $\h$ is easily calculated as:
$$
 \bar\m\indb(\f\in\O_L:\h_1(\f)=\h_1,\ldots,\h_L(\f)=\h_L)
  =\indfn{|\h_1|\leq M}\frac{\nep{-\b\sum_{i=2}^L|\h_i|-W\indc(\b,T(\h))}}{\bar Z\indb(\b)}.
$$
Observe that $\h(\f)$ defined
in \equ(v.defeta) is the vector of the discrete derivatives of the configuration $\f\in\O_L$.
The study of the SOS interface can be carried out using the variable $\f$ 
or $\h$ indifferently.
We will use the last one in the sequel.
Let $\bnu$ be the distribution of $\h$, \ie
$$
 \bnu(\h)
 \equiv\indfn{|\h_1|\leq M}\frac{\nep{-\b\sum_{i=2}^L|\h_i|-W\indc(\b,T(\h))}}
 {\bar Z\indb(\b)}.
\Eq(v.defn)
$$ 
The expected value operator with respect to $\bnu$ will be denoted with $\omean(\cdot)$,
while the {\it covariance\/} form will be denoted by $\omean(\cdot,\cdot)$ \ie
$$
 \omean(f,g)\equiv\omean\pquad{\ptond{f-\omean(f)}\ptond{g-\omean(g)}},
$$
where $f,g\in L^2(\bnu)$.
On the same Hilbert space is defined the quadratic form:
$$
 \oce(f,f)\equiv
  \omean\pquad{\indfn{\h_1<M}(\de_{2,1}f)^2}+\sum_{k=2}^{L-1}\omean\pquad{(\de_{k+1,k}f)^2}+\omean\pquad{(\de^+_Lf)^2},
$$
where
$$
 \eqalign{
  (\de_h^+ f)(\h)&=f(\h+\d_h)-f(\h)\cr
  (\de_h^- f)(\h)&=f(\h-\d_h)-f(\h)\qquad\qquad\qquad\qquad\qquad h,k=1,\ldots,L\cr
  (\de_{h,k}f)(\h)&=f(\h+\d_h-\d_k)-f(\h).\cr
 }
$$
We will use this form to estimate the spectral gap of $\bar G\indb$:

\nproclaim Lemma [minmax1].
There exists two constants $K_1(\b)$ and $K_2(\b)$ such that
$$
 K_1\inf_{f\in L^2(\bnu)}{\oce(f,f)\over\omean(f,f)}
  \leq\l_1(G\indb)
  \leq K_2\inf_{f\in L^2(\bnu)}{\oce(f,f)\over\omean(f,f)}.
\Eq(g-e0)
$$

\Pro
Using the fact that $\bar G\indb$ is self adjoint in  $L^2(\bar\m\indb)$ 
it simple to check that
$$
 \bar\cG\indb(f,f)\equiv -\scalprod{\bar G\indb f}{f}{L^2(\bar\m\indb)}=
  \ov2\sum_{k=1}^{L}\sum_{\f\in\O_L}\bar\m\indb(\f)\bar c\indb(\f,\f+\d_k)\pquad{(\de_k^+f)(\f)}^2.
$$
If we recall the definition of the jump rates \equ(defrates3) a simple 
calculation shows that there 
exists $C_1(\b)$ and $C_2(\b)$ such that
$$
 C_1 \bar\m\indb(\f)\indfn{\f_1<M}\leq \bar\m\indb(\f)\bar c\indb(\f,\f+\d_1)\leq C_2 \bar\m\indb(\f)\indfn{\f_1<M}
$$
and
$$
 C_1 \bar\m\indb(\f)\leq \bar\m\indb(\f)\bar c\indb(\f,\f+\d_k)\leq C_2 \bar\m\indb(\f)
        \qquad\qquad
        k=2,\ldots,L,
$$
for every $\f\in\O_L$ and $L,M>0$.
This means that $\bar\cG\indb(f,f)$ can be bounded from above an from below by
$$
 \sum_{\f\in\O_L}\bar\m\indb(\f)\pgraf{\pquad{\indfn{\f_1<M}(\de_1^+f)(\f)}^2+\sum_{k=1}^{L}\pquad{(\de_k^+f)(\f)}^2}
$$
multiplyed by $C_1/2$ and $C_2/2$ respectively.
Now we use the change of variable $\f=T\h$ to obtain:
$$
 \sum_{\f\in\O_L}\bar\m\indb(\f)\pgraf{\pquad{\indfn{\f_1<M}(\de_1^+f)(\f)}^2+\sum_{k=1}^{L}\pquad{(\de_k^+f)(\f)}^2}
 =\oce(f^*,f^*)
$$
where $f^*(\h)\equiv f(T\h)$.
This implies
$$
 {C_1\over 2}\inf_{f^*\in L^2(\bnu)}
  {\oce(f^*,f^*)\over\omean(f^*,f^*)}
   \leq\inf_{f\in L^2(\bar\m\indb)}
    {\bar\cG\indb(f,f)\over\bar\m\indb(f,f)}
     \leq {C_2\over 2}\inf_{f^*\in L^2(\bnu)}
    {\oce(f^*,f^*)\over\omean(f^*,f^*)},
\Eq(r:1)
$$
where $\bar\m\indb(f,f)\equiv\int d\bar\m\indb(\f)f^2(\f)-\pquad{\int d\bar\m\indb(\f)f(\f)}^2=\omean(f^*,f^*)$.
Observing that
$$
  \inf_{g\in L^2(\bnu)}
  {\oce(g,g)\over\omean(g,g)}=
  \inf_{g\in L^2(\bnu)\atop g\perp{\bf 1}}
  {\oce(g,g)\over\omean(g^2)},
        \qquad\qquad
        \inf_{g\in L^2(\bar\m\indb)}
    {\bar\cG\indb(g,g)\over\bar\m\indb(g,g)}=
  \inf_{g\in L^2(\bar\m\indb)\atop g\perp{\bf 1}}
    {\bar\cG\indb(g,g)\over\int d\bar\m\indb(\f)g^2(\f)},
$$
by the variational characterization of the spectral gap and \equ(r:1) we have \equ(g-e0).
\qed

We can use this lemma to prove Proposition~\thm[v5].
In fact by \equ(g-e0)
it is easy to show that \equ(v5.1) is 
equivalent to the existence of $C_1(\b)>0$ such that the Poincar\'e inequality
$$
 \omean(f,f)\leq C_3(\b)(L\lor M^2)L\oce(f,f)
\Eq(v.poi)
$$
holds for every $f\in L^2(\bnu)$.
The key step of the proof of this inequality is contained in the following result.

\nproclaim Proposition [gla].
There exists  $\bar\b>0$ such that for every $\b\geq\bar\b$ it is 
possible to find
$K_1(\b)>0$ such that
$$
 \omean(f,f|\h_1)\leq K_1\sum_{k=2}^L\omean\pquad{(\de_k^+f)^2|\h_1}
\Eq(gla.1)
$$
for any $f\in L^2(\bnu)$, $L>0$ and $\h_1\in[-M,M]\cap\integer$.

\noindent
This proposition shows the perturbative approach of the proof.
In fact if $W\indc=0$ the measure $\bnu(\cdot|\h_1)$ is a 
product measure.
In this case, Proposition~\thm[gla] says that there exists a uniformly positive spectral gap
for a random walk in $\integer^{L-1}$ in which each component of the walk is independent 
from the others.
It is well known that this gap exists if each component exhibits by itself 
a uniformly positive spectral gap, and the existence of this one site 
spectral gap is easily proved.
Because for $\b\to+\infty$ the perturbation $W\indc$ goes to $0$, the result should be 
true also for large values of $\b$.

Before proving the key result Proposition~\thm[gla] we want to show how, 
from this result, follows \equ(v.poi).

\Prot{of Proposition \thm[v5]}
Fix $f\in L^2(\bnu)$, a simple calculation yields:
$$
 \omean(f,f)=\omean\pquad{\omean(f,f|\h_1)}+\omean\pquad{\omean(f|\h_1),\omean(f|\h_1)}.
\Eq(v.poi.1)
$$
The first term on the right hand side of this equation can be bounded using
Proposition~\thm[gla].
For the second term notice that $\h_1$ is uniformly distributed in
$[-M,M]\cap\integer$.
It is well known (see [1] for example) that this implies that there 
exists $C_1(\b)>0$ such that the 
one-site Poincar\'e inequality
$$
 \omean\ptond{g,g}\leq C_1M^2\omean\pquad{\indfn{\h_1<M}(\de_1^+g)^2}
$$
holds for every $g=g(\h_1)\in L^2(\bnu)$ and $M>0$.
In particular this formula is true for $g(\h_1)=\omean(f|\h_1)$.
Because the measure
$\bnu(\cdot|\h_1)$ does not depend on $\h_1$ if $-M\leq \h_1< M$, is simple 
to check that
$$
 \de_1^+\omean(f|\h_1)=\omean(\de_1^+f|\h_1).
$$
In conclusion we obtain
$$
 \omean\pquad{\omean(f|\h_1),\omean(f|\h_1)}
  \leq C_1M^2\omean\pquad{\indfn{\h_1<M}(\de_1^+f)^2}.
$$
We can use this and \equ(gla.1) in the left hand side of
\equ(v.poi.1) to get:
$$
 \omean(f,f)\leq C_2(\b)\pgraf{\sum_{k=2}^L\omean\pquad{(\de_k^+f)^2}
  +M^2\omean\pquad{\indfn{\h_1<M}(\de_1^+f)^2} }.
\Eq(v.poi.2)
$$
Now notice that
$$
 (\de_k^+f)(\h)=\sum_{i=k}^{L-1}(\de_{i+1,i}f)(\h+\d_{i+1})+(\de_L^+f)(\h)
$$
and that the change of coordinates $\h\mapsto\h+\d_j$ has bounded Jacobian.
Thus \equ(v.poi.2) gives:
$$
 \omean(f,f)\leq C_3(\b)(M^2\lor L)L\pgraf{\omean\pquad{\indfn{\h_1<M}(\de_{2,1}f)^2}+\sum_{k=2}^L\omean\pquad{(\de_{k+1,k}f)^2}
  +\omean\pquad{(\de_L^+f)^2} },
$$
which is the same as \equ(v.poi).
\qed

The remaining part of this section is devoted to the proof of Proposition~\thm[gla].
It is convenient to introduce some extra notation.
Recall that $\h_1$ is  independent from $\h_2,\ldots,\h_L$ and that 
$W_L^{\infty}(\b,T(\h))$ does not depend on $\h_1$.
This implies that
$$
 \bnu(\h|\h_1)
  =\frac{(2M+1)\nep{-\b\sum_{i=2}^L|\h_i|-W_L^{\infty}(\b,T(\h))}}
   {\bar Z\indb(\b)},
$$
for every $|\h_1|\leq M$.
If we define $\h_1^\prime\equiv\h_2$,
$\h_2^\prime\equiv\h_3$,\dots, $\h_{L-1}^\prime\equiv\h_L$ and
$\hat W_{L-1}(\h^\prime)\equiv W_L^\infty(T(\h))$,
the last expression becomes:
$$
 \bnu(\h|\h_1)
  =\frac{(2M+1)\nep{-\b\sum_{i=1}^{L-1}|\h_i^\prime|-\hat W_{L-1}(\b,\h^\prime)}}
   {\bar Z\indb(\b)}. 
$$
Define on $\O_L$ the probability measure
$$
 \bnl(\h)\equiv\frac{\nep{-\b\sum_{i=1}^{L}|\h_i|-\hat W_{L}(\b,\h)}}
   {\bar Z_L(\b)},
 \qquad\qquad
 \bar Z_L(\b)\equiv\sum_\h\nep{-\b\sum_{i=1}^{L}|\h_i|-\hat W_{L}(\b,\h)} 
$$
and denote with $\oml(\cdot)$ the expected value with respect to this measure.
We restate Proposition~\thm[gla] as

\nproclaim Proposition [gla1].
There exists  $\bar\b>0$ such that for every $\b\geq\bar\b$ it is 
possible to find
$K_1(\b)>0$ such that
$$
 \oml(f,f)\leq K_1(\b)\sum_{k=1}^L\oml\pquad{(\de_k^+f)^2},
\Eq(goalgla1)
$$
for any $f\in L^2(\bnl)$ and $L>0$.

\noindent
We  will prove this result using the {\it martingale approach\/} outlined in [5].

Define the subsets $\a_j$ of $\natural$ in the following way
$$
 \a_j
   \equiv\cases{
     \pgraf{j,\ldots,L}& if $1\leq j\leq L$,\cr
     \emptyset& if $j=L+1$.\cr
    }
$$
For every $j=1,\ldots,L$ the restriction of $\h$ to the set $\a_j$ is 
denoted by $\h_{\a_{j}}\equiv(\h_j,\ldots,\h_L)$.
If we define $f_j\equiv\oml(f|\h_{\a_{j}})$
it is simple to check that:
$$
 \oml(f,f)=\sum_{i=1}^L\oml\pquad{\oml(f_j,f_j|\h_{\a_{j+1}})}.
\Eq(vv.1)
$$
On the right hand side of this formula there is a sum of expected value of 
conditional variances.
Notice that the random variable $f_j$, by definition, depends only on $\h_{\a_{j}}$.
So if $\h_{\a_{j+1}}$ is fixed, it depends only on $\h_j$ (we will say 
that $f$ {\it is local\/} in $j$).
This means that each of the variances on the right hand side of 
\equ(vv.1) is the variance of a local function.

The method we will use to prove \equ(goalgla1) consists of two steps
The first step is to show that the marginal in $\h_j$ of the measure 
$\bnl(\cdot|\h_{\a_{j+1}})$ exhibits a positive spectral gap uniformly in $L>0$.

\nproclaim Lemma [spg].
There exists $\bar\b>0$ such that for every $\b>0$ there exists $K_1(\b)>0$ 
so that
$$
 \oml(f,f|\h_{\a_{j}})\leq K_1(\b)\oml\pquad{(\de_j^+f)^2|\h_{\a_{j}}},
\Eq(vax1)
$$
for every $L>0$ and $f\in L^2(\bnl(\cdot|\h_{\a_{j}}))$ local in $j$.

Because of this lemma \equ(vv.1) becomes
$$
 \oml(f,f)\leq K_1(\b)\sum_{i=1}^L\oml\pquad{(\de_j^+f_j)^2}.
\Eq(stp)
$$
The second step is to show that the right hand side of \equ(stp) can be bounded from 
above by the correct quadratic form:

\nproclaim Lemma [v1].
There exists $\bar\b$ such that
$$
 \sum_{i=1}^L\oml\pquad{(\de_i^+f_i)^2}\leq 4\sum_{i=1}^L\oml\pquad{(\de_i^+f)^2},
\Eq(v1.1)
$$
for every $\b\geq\bar\b$, $L>0$ and $f\in L^2(\bar\n)$.

\noindent
If we use Lemma~\thm[v1] in \equ(stp) we obtain immediately \equ(goalgla1).

In order to prove Lemma~\thm[spg] and Lemma~\thm[v1] we need a preliminary result.

\nproclaim Lemma [v3].
There exists $\bar\b>0$ such  that for every $\b>\bar\b$
$$
 \exp\pquad{-\b\sign(\h_j)-8\nep{-m(\b)}}
  \leq\frac{\bnl(\h_j+\d_j|\h_{\a_{j+1}})}{\bnl(\h_j|\h_{\a_{j+1}})}
   \leq\exp\pquad{-\b\sign(\h_j)+8\nep{-m(\b)}}
\Eq(v3.1)
$$
and
$$
 \exp\pquad{-16\nep{-m(\b)(i-j)}}
  \leq\frac{\bnl(\h_j|\h_{\a_{j+1}}+\d_i)}{\bnl(\h_j|\h_{\a_{j+1}})}
   \leq\exp\pquad{16\nep{-m(\b)(i-j)}},
\Eq(v3.1.1)
$$
for every $j=1,\ldots,L-1$, $\h\in\O_L$.

\Pro
The proof is divided in several steps for purpose of clarity.
To keep notation simple we will write $\D(\h)$ instead of $\D(T\h)$
and
$$
 \cS(\h)\equiv\pgraf{\L\sset V_L^\infty:\L\cap\D\ptond{\h}\not=\emp}.
$$
Recall that
$$
 \hat W_{L}(\b,\h)=\sum_{\L\in\cS(\h)}\F(\b,\L).
$$
Define for $z\in\real$ the line 
$r(z)\equiv\pgraf{(x,y)\in\real^2:x=z}$ and for fixed $j=1,\ldots,L-1$
$$
 \cS_j(\h)\equiv\cS(\h)\setminus\pgraf{\L\in\cS(\h):\L>r(j)}.
$$

\proclaim Step 1.
Define
$$
 \hat W_j(\h)\equiv\sum_{\L\in\cS_j(\h)}\F(\b,\L).
$$
Then
$$
 \bnl(\h_j|\h_{\a_{j+1}})
  =\frac{\nep{-\b|\h_j|}\sum_{\h_1,\ldots,\h_{j-1}}\nep{-\b\sum_{k=1}^{j-1}|\h_k|-\hat W_j(\h)}}
   {\sum_{\h_1,\ldots,\h_{j}}\nep{-\b\sum_{k=1}^{j}|\h_k|-\hat W_j(\h)}}.
\Eq(v3.2)
$$

\Prot{of Step 1}
Define
$$
 \hat W_j^c(\h)\equiv\hat W_L(\h)-\hat W_j(\h)=\sum_{\L\in\cS(\h)\setminus\cS_j(\h)}\F(\b,\L),
$$
we claim that $\hat W_j^c$ does not depend on $\h_1,\ldots,\h_j$.
If we assume this we obtain:
$$
 \eqalign{
  &\quad \bnl(\h_j|\h_{\a_{j+1}})
   =\frac{\nep{-\b|\h_j|}\sum_{\h_1,\ldots,\h_{j-1}}\nep{-\b\sum_{k=1}^{j-1}|\h_k|-\hat W_L(\h)}}
    {\sum_{\h_1,\ldots,\h_{j}}\nep{-\b\sum_{k=1}^{j}|\h_k|-\hat W_j(\h)}}=\cr
  &=\frac{\nep{-\b|\h_j|}\sum_{\h_1,\ldots,\h_{j-1}}\nep{-\b\sum_{k=1}^{j-1}|\h_k|-\hat W_j(\h)}}
   {\sum_{\h_1,\ldots,\h_{j}}\nep{-\b\sum_{k=1}^{j}|\h_k|-\hat W_j(\h)}},\cr
 }
$$
\ie \equ(v3.2).

In order to prove that $\hat W_j^c(\h)$ does not depend on $\h_1,\ldots,\h_j$
it suffices to show that
$$
 \hat W_j^c(\h+h\d_k)=\hat W_j^c(\h)
$$
for every $h\in\integer$ and $k=1,\ldots,j$.
Define
$$
 T_h:\L\ni\cS(\h)\setminus\cS_j(\h)\mapsto\L+h\in\cS(\h+h\d_k)\setminus\cS_j(\h+h\d_k).
$$
This map is bijective for every 
$k=1,\ldots,j$. Furthermore because of the translation invariance \equ(dec2) of $\F(\b,\cdot)$
we have:
$$
 \eqalign{
  &\quad \hat W_j^c(\h)
   =\sum_{\L\in\cS(\h)\setminus\cS_j(\h)}\F(\b,\L)
   =\sum_{\L\in\cS(\h)\setminus\cS_j(\h)}\F(\b,T_h\L)=\cr
  &=\sum_{\L\in\cS(\h+h\d_k)\setminus\cS_j(\h+h\d_k)}\F(\b,\L)
   = \hat W_j^c(\h+h\d_k).\cr
 }
$$
The \equ(v3.2) is proved.
\qed

From Step~1 we obtain
$$
 \frac{\bnl(\h_j+\d_j|\h_{\a_{j+1}})}{\bnl(\h_j|\h_{\a_{j+1}})}
  =\nep{-\b(|\h_j+1|-|\h_j|)}
   \frac{\sum_{\h_1,\ldots,\h_{j-1}}\nep{-\b\sum_{k=1}^{j-1}|\h_k|-\hat W_j(\h+\d_j)}}
   {\sum_{\h_1,\ldots,\h_{j-1}}\nep{-\b\sum_{k=1}^{j-1}|\h_k|-\hat W_j(\h)}}
$$
and:
$$
 \frac{\bnl(\h_j|\h_{\a_{j+1}}+\d_i)}{\bnl(\h_j|\h_{\a_{j+1}})}
  =\frac{\sum_{\h_2,\ldots,\h_{j-1}}\nep{-\b\sum_{i=2}^{j-1}|\h_i|-\hat W_j(\h+\d_i)}}
   {\sum_{\h_2,\ldots,\h_{j-1}}\nep{-\b\sum_{i=2}^{j-1}|\h_i|-\hat W_j(\h)}}
    \frac{\sum_{\h_2,\ldots,\h_{j}}\nep{-\b\sum_{i=2}^{j}|\h_i|-\hat W_j(\h)}}
    {\sum_{\h_2,\ldots,\h_{j}}\nep{-\b\sum_{i=2}^{j}|\h_i|-\hat W_j(\h+\d_i)}},
$$
for every $i>j$.
By these inequalities we obtain that to prove \equ(v3.1) and \equ(v3.1.1) we have only to show that
$$
 \sup_\h|\hat W_j(\h+\d_i)-\hat W_j(\h)|\leq 8\nep{-m(\b)(i+1-j)}
\Eq(v3.goal)
$$
for every $i\geq j$.

\proclaim Step 2.
Define
$$
 \D_i(\h)\equiv\D(\h)\cap\pgraf{(x,y)\in\real^2: x>i}
$$
and
$$
 \cS_{i,j}(\h)\equiv\pgraf{\L\sset V_L^\infty:\L\cap\D_i(\h)\not=\emp,\ \L\cap r(j)\not=\emp}.
$$
Then:
$$
 \hat W_j(\h+\d_i)-\hat W_j(\h)
  =\sum_{\L\in\cS_{i,j}(\h+\d_i)}\F(\b,\L)-\sum_{\L\in\cS_{i,j}(\h)}\F(\b,\L).
\Eq(v3.3)
$$

\Prot{of Step 2}
Let $\L\in\cS_j(\h)$ be such that it intersects $\D(\h)$ on the left of $i$.
Then it also intersects $\cS_j(\h+\d_i)$ in the same points.
On the contrary if $\L\in\cS_j(\h)$ intersects $\D(\h+\d_i)$ on the left of $i$ then it intersects $\cS_j(\h)$.
In conclusion
$$
 \pgraf{\L\in\cS_j(\h):
  \vcenter{
   \hbox{$\L$ intersects $\D(\h)$}
   \hbox{on the left of $i$}
  }
 }
  =\pgraf{\L\in\cS_j(\h+\d_i):
  \vcenter{
   \hbox{$\L$ intersects $\D(\h+\d_i)$}
   \hbox{on the left of $i$}
  }
 }.
$$
By this relation we obtain that we can clear from the difference
$$
 \hat W_j(\h+\d_i)-\hat W_j(\h)
  =\sum_{\L\in\cS_j(\h)}\F(\b,\L)-\sum_{\L\in\cS_j(\h+\d_i)}\F(\b,\L),
$$
all the terms $\F(\b,\L)$ such that $\L$ intersects $\D(\h)$,
or $\D(\h+\d_i)$, on the left of $i$.
It follows that the sums are actually only on the $\L$ which neither 
intersects $\D(\h)$ on the left of $i$ nor intersects
$\D(\h+\d_i)$ on the left of $i$.
Because these $\L$ have to intersect
$\D(\h)$, the intersection is on the right of $i$.
This proves \equ(v3.3).
\qed

For any $S\sset\dual$
we will say that $p\in S$ is {\it $+$unstable\/} if:
$$
 p\notin (S+{\bf e}_y),
$$
where ${\bf e}_y=(0,1)\in\real^2$.
Similarly we will say that
$p\in S$ is {\it $-$unstable} if:
$$
 p\notin (S-{\bf e}_y).
$$
The classes of points $+$unstable and $-$unstable of the set $S$
will be denoted respectively with $\cI^+(S)$ and $\cI^+(S)$.

\proclaim Step 3.
Define
$$
 \eqalign{
  \cS_{i,j}^+(\h)&\equiv\cS_{i,j}(\h)\cap\cI^+(\D_i(\h))\cr
  \cS_{i,j}^-(\h)&\equiv\cS_{i,j}(\h)\cap\cI^-(\D_i(\h)).\cr
 }
$$
Then
$$
 \hat W_j(\h+\d_i)-\hat W_j(\h)
  =\sum_{\L\in\cS_{i,j}^-(\h+\d_i)}\F(\b,\L)
   -\sum_{\L\in\cS_{i,j}^+(\h)}\F(\b,\L).
\Eq(v3.4)
$$

\Pro
It simple to check that:
$$
 \L\in\cS_{i,j}(\h+\d_i),\quad\L\notin\cI^-(\D_i(\h+\d_i))
 \qquad\Longleftrightarrow\qquad
 \L\in\cS_{i,j}(\h),\quad\L\notin\cI^+(\D_i(\h)).
$$
This implies $\cS_{i,j}(\h+\d_i)\setminus\cI^-(\D_i(\h+\d_i))=\cS_{i,j}(\h)\setminus\cI^+(\D_i(\h))$,
which proves \equ(v3.4).
\qed

\proclaim Step 4.
For any $\h\in\O_L$ and $k=1/2,1+1/2,2+1/2,\ldots$
$$
 \Norm{\cI^\pm(\D(\h))\cap r(k)}\leq 2.
$$

\Pro
If $\D=\D(\h)$ it is clear that $\D=\overline\D\cup\underline\D$ where
$$
 \overline\D\equiv\pgraf{p\in\D:\hbox{$p$ is above $\G(T\h)$}}
  \qquad\qquad\qquad
  \underline\D\equiv\pgraf{p\in\D:\hbox{$p$ is below $\G(T\h)$}}.
$$
Notice that in general
$\cI^\pm(A\cup B)\sset\cI^\pm(A)\cup\cI^\pm(B)$.
Thus to prove
\equ(v3.4) we have only to show that:
$$
 \Norm{\cI^\pm(\overline\D)}=1\qquad\qquad\hbox{and}\qquad\qquad\Norm{\cI^\pm(\underline\D)}=1.
$$
This fact can be easily checked by using geometric considerations.
\qed

We are finally in a position to prove \equ(v3.goal).
Notice that if $\L\in\cS_{i,j}^+$ then:
\item{$i)$} $\L$ contains $p\in\cI^+(\D_i(\h))$;
\item{$ii)$} $\L$ intersects $r(j)$, thus because $\L\sset\dual$, intersects $r(j-1/2)$.

\noindent
It follows that:
$$
 \eqalign{
 &\quad\sum_{\L\in\cS_{i,j}^+(\h)}\F(\b,\L)
  \leq\sum_{k=i}^{+\infty}
   \sum_{p\in\cI^+(\D_i(\h))\cap r(k+1/2)}
    \sum_{\st\L\ni p\atop\st\L\cap r(j-1/2)\not=\emp}\F(\b,\L)\leq\cr
 &\leq2\sum_{k=i}^{+\infty}\sum_{\st\L\ni p\atop\st\diam{\L}\geq k-j+1}\F(\b,\L)\leq
  2\sum_{k=i}^{+\infty}\nep{-m(\b)(k-j+1)}
  \leq 4\nep{-m(\b)(k-j+1)}.
}
$$
Where  $\b$ is large enough and \equ(100) has been used.
From this estimate we obtain \equ(v3.goal) that, as we noticed before, implies
\equ(v3.1) and \equ(v3.1.1).
\qed

Lemma~\thm[v3] can be used to prove the one-site spectral gap Lemma~\thm[spg].
In fact \equ(v3.1) shows that for $\b$ large enough (recall that 
$m(\b)\to+\infty$ for $\b\to+\infty$) the measure 
$\bnl(\cdot|\h_{\a_j})$ exhibits an ``inward drift'' (see [11]).

\Prot{of Lemma~\thm[v3]}
It follows immediately from [3] and [10].
\qed

Now we turn to the proof of Lemma~\thm[v1].
We need a technical result

\nproclaim Lemma [v2].
For every $i=2,\ldots, L$ we have:
$$
 \de_i^+f_i=\oml(\de_if|\h_{\a_i})+\sum_{j=1}^{i-1}\oml(f^+_j,V_{i,j}|\h_{\a_i}),
\Eq(v2.1)
$$
where $f^+_j(\h)\equiv f_j(\h+\d_i)$ and
$V_{i,j}(\h)\equiv{\bnl(\h_j|\h_{\a_{j+1}}+\d_i)\over \bnl(\h_j|\h_{\a_{j+1}})}-1$.

\Pro
For every $k=1,\ldots,L$ define
$$
 F_k(\h)\equiv\bnl(\h_k|\h_{\a_{k+1}}).
$$
$F_k$ is the marginal in $\h_k$ of $\bnl(\cdot|\h_{\a_{k+1}})$.
Notice that
$$
 f_i=\oml(f|\h_{\a_i})=\oml(\oml(\cdots\oml(f|\h_{\a_2})\cdots|\h_{\a_{i-1}})|\h_{\a_i})
$$
for any $i=2,\ldots, L$.
Thus
$$
 f_i=\sum_{\h_1,\ldots,\h_{i-1}}f\prod_{k=1}^{i-1}F_k.
\Eq(v2.2.1)
$$
To compute $\de_i^+f_i$ we have to calculate the (discrete) derivative of 
a product.
We will use the following formula for the derivative of products:
$$
 (\de_i^+gh)=(\de_i^+g)h^++g(\de_i^+h),
\Eq(v2.3)
$$
where $g^+(\h)\equiv g(\h+\d_i)$.
Using \equ(v2.3) we obtain:
$$
 \de_i^+\Big(f\prod_{k=1}^{i-1}F_k\Big)=(\de_i^+f)\prod_{k=1}^{i-1}F_k
  +\sum_{j=1}^{i-1}f^+(\de_i^+F_j)\Big(\prod_{l=1}^{j-1}F_l^+\Big)
   \Big(\prod_{k=j+1}^{i-1}F_k\Big).
$$
This relation and \equ(v2.2.1) give:
$$
 \eqalign{
  \de_i^+f_i
  &=\sum_{\h_1,\ldots,\h_{i-1}}(\de_i^+f)\prod_{k=1}^{i-1}F_k
   +\sum_{\h_1,\ldots,\h_{i-1}}f^+\sum_{j=1}^{i-1}(\de_i^+F_j)\Big(\prod_{l=1}^{j-1}F_l^+\Big)
   \Big(\prod_{k=j+1}^{i-1}F_k\Big)=\cr
  &=\oml(\de_i^+f|\h_{\a_i})
   +\sum_{j=1}^{i-1}\sum_{\h_j,\ldots,\h_{i-1}}(\de_i^+F_j)\Big(\prod_{k=j+1}^{i-1}F_k\Big)
    \sum_{\h_1,\ldots,\h_{j-1}}f^+\Big(\prod_{l=1}^{j-1}F_l^+\Big)=\cr
  &=\oml(\de_i^+f|\h_{\a_i})
   +\sum_{j=1}^{i-1}\sum_{\h_j,\ldots,\h_{i-1}}f^+_j(\de_i^+F_j)\Big(\prod_{k=j+1}^{i-1}F_k\Big)=\cr
  &=\oml(\de_i^+f|\h_{\a_i})
   +\sum_{j=1}^{i-1}\sum_{\h_j,\ldots,\h_{i-1}}
    f^+_j\ptond{{F_j^+\over F_j}-1}\Big(\prod_{k=j}^{i-1}F_k\Big)=\cr
  &=\oml(\de_i^+f|\h_{\a_i})+\oml(f^+V_{i,j}|\h_{\a_i}).\cr
 }
$$
To prove \equ(v2.1) it remains to observe that by definition $\oml(V_{i,j}|\h_{\a_i})=0$.
\qed

\Prot{of Lemma~\thm[v1]}
We borrow the basic idea of the proof from [9].
We will show that for $\b$ large enough
$$
 \sum_{i=1}^L\oml\pquad{(\de_i^+f_i)^2}\leq 2\sum_{i=1}^L\oml\pquad{(\de_i^+f)^2}
  +\ov2\sum_{i=1}^L\oml\pquad{(\de_i^+f_i)^2}, 
\Eq(goalv1)
$$
for every $f\in L^2(\bnl)$.

Fix $f\in L^2(\bnl)$ and $i>1$.
By Lemma~\thf[v2]
$$
 (\de_i^+f_i)^2
  \leq 2\oml\pquad{(\de_i^+f)^2|\h_{\a_i}}
   +2\Big[\sum_{j=1}^{i-1}\oml(f_j^+,V_{i,j}|\h_{\a_i})\Big]^2.
\Eq(v1.2)
$$
We have to estimate the second term on the right hand side of this relation.
By Schwartz inequality and Lemma~\thf[v3]
we obtain
$$
 \oml(f_j^+,V_{i,j}|\h_{\a_i})\leq2\e^{i-j}\oml(f_j^+,f_j^+|\h_{\a_i})^{\ov2},
$$
where $\e=\e(\b)\to 0$ for $\b\to+\infty$.
Thus if $\b$ is large enough
$$
 \eqalign{
 &\quad\Big[\sum_{j=1}^{i-1}\oml(f_j^+,V_{i,j}|\h_{\a_i})\Big]^2
  \leq 4\Big[\sum_{j=1}^{i-1}\e^{i-j}\oml(f_j^+,f_j^+|\h_{\a_i})^{\ov2}\Big]^2\leq\cr
 &\leq 4\Big(\sum_{j=1}^{i-1}\e^{i-j}\Big)\sum_{j=1}^{i-1}\e^{i-j}\oml(f_j^+,f_j^+|\h_{\a_i})
  \leq 8\e\sum_{j=1}^{i-1}\e^{i-j}\oml(f_j^+,f_j^+|\h_{\a_i}).\cr
 }
\Eq(v1.3)
$$
It simple to check that
$$
 \oml(f_j^+,f_j^+|\h_{\a_i})=\sum_{s=j}^{i-1}\oml\pquad{\oml(f_s^+,f_s^+|\h_{\a_{s+1}})|\h_{\a_i}}.
$$
So by \equ(vax1) we know that, for $\b$ large enough, there 
exists
$C_1>0$ so that:
$$
 \eqalign{
  &\quad\oml(f_j^+,f_j^+|\h_{\a_i})
   \leq C_1
    \sum_{s=j}^{i-1}\oml\pquad{\oml\ptond{(\de_s^+f_s^+)^2|\h_{\a_{s+1}}}|\h_{\a_i}}=\cr
  &=C_1\sum_{s=j}^{i-1}\oml\pquad{(\de_s^+f_s^+)^2|\h_{\a_i}}.\cr
 } 
$$
By using this bound in \equ(v1.3) we get
$$
 \Big[\sum_{j=1}^{i-1}\oml(f_j^+,V_{i,j}|\h_{\a_i})\Big]^2
  \leq \e C_2\sum_{j=1}^{i-1}\e^{i-j}\sum_{s=j}^{i-1}\oml\pquad{(\de_s^+f_s^+)^2|\h_{\a_i}}.
\Eq(v1.4)
$$
Exchanging the sums on the right hand side of \equ(v1.4) we have
$$
 \eqalign{
  &\quad\sum_{j=1}^{i-1}\e^{i-j}\sum_{s=j}^{i-1}\oml\pquad{(\de_s^+f_s^+)^2|\h_{\a_i}}
   =\sum_{s=1}^{i-1}\oml\pquad{(\de_s^+f_s^+)^2|\h_{\a_i}}\sum_{j=1}^{s}\e^{i-j}\leq\cr
  &\leq 2\sum_{s=1}^{i-1}\e^{i-s}\oml\pquad{(\de_s^+f_s^+)^2|\h_{\a_i}}.\cr
 }
$$
This implies
$$
 \Big[\sum_{j=1}^{i-1}\oml(f_j^+,V_{i,j}|\h_{\a_i})\Big]^2
  \leq \e C_3\sum_{s=1}^{i-1}\e^{i-s}\oml\pquad{(\de_s^+f_s^+)^2|\h_{\a_i}}.
$$
and by \equ(v1.2):
$$
 (\de_i^+f_i)^2
  \leq 2\oml\pquad{(\de_i^+f)^2|\h_{\a_i}}
   +\e C_4\sum_{s=1}^{i-1}\e^{i-s}\oml\pquad{(\de_s^+f_s^+)^2|\h_{\a_i}}.
$$
Taking expected value on both sides of this relation and recalling that 
the change of variable
$\h\mapsto\h-\d_i$ has a bounded Jacobian (see Lemma~\thm[v3]) 
we obtain
$$
 \oml\pquad{(\de_i^+f_i)^2}
  \leq 2\oml\pquad{(\de_i^+f)^2}
   +\e C_4\sum_{s=1}^{i-1}\e^{i-s}\oml\pquad{(\de_s^+f_s)^2}.
$$
We sum this relation for $i=2,\ldots,L$.
An elementary computation gives
$$
 \eqalign{
  \sum_{i=2}^L\oml\pquad{(\de_i^+f_i)^2}
  &\leq 2\sum_{i=2}^L\oml\pquad{(\de_i^+f)^2}
   +C_4\e\sum_{i=2}^L\sum_{s=1}^{i-1}\e^{i-s}\oml\pquad{(\de_s^+f_s)^2}=\cr
  &=2\sum_{i=2}^L\oml\pquad{(\de_i^+f)^2}
   +C_4\e\sum_{s=1}^{L-1}\sum_{i=s+1}^{L}\e^{i-s}\oml\pquad{(\de_s^+f_s)^2}\leq\cr
  &\leq 2\sum_{i=2}^L\oml\pquad{(\de_i^+f)^2}
   +2C_4\e^2\sum_{s=1}^{L-1}\oml\pquad{(\de_s^+f_s)^2}.\cr
 }
$$
Recalling that $\de_i^+f_1=\de_i^+f=\oml(\de_i^+f|\h_{\a_1})$, this implies:
$$
 \sum_{i=1}^L\oml\pquad{(\de_i^+f_i)^2}
  \leq 2\sum_{i=1}^L\oml\pquad{(\de_i^+f)^2}
   +\e^2 C_5\sum_{s=1}^{L}\oml\pquad{(\de_s^+f_s)^2}.
$$
To conclude the proof of \equ(goalv1)
we choose $\bar\b>0$ such that 
$\b\geq\bar\b$ implies $C_5\e^2(\b)<1/2$.
\qed

\fine


\expandafter
 \ifx\csname sezioniseparate\endcsname\relax
  \fi

\numsec=-1
\numfor=1
\numtheo=1
\pgn=1

\beginsection A. Appendix

\nproclaim Lemma [v8].
Define $A\equiv\pgraf{\f\in\O_L:\normm{\f}_{\infty}\leq (1-\e)L/2}$
and $B\equiv\pgraf{\f\in\O_L:\normm{\f}_{\infty}\leq\a L}$
where $\e\in\ptond{0,1/100}$ and $\a\in\ptond{0,1/4}$.
Then there exists $\bar\b>0$ such that for every $\b>\bar\b$ there 
exists  $K_1(\b)>0$ so that
$$
 \sup_{\f\in A}{\m\inda(\f|B)\over\bar\m\inda(\f)}\leq \a^{-1}K_1.
\Eq(v8.1)
$$

\Pro
An elementary calculation shows that:
$$
 {\m\inda(\f|B)\over\bar\m\inda(\f)}
 \leq\nep{\sup_{\f\in B}|W\indc(\f)-W\inda(\f)|}\bar\m\inda(B)^{-1},
\Eq(v8.2)
$$
for every $\f\in B$.
Observe that:
$$
 W\indc(\b,\f)-W\inda(\b,\f)=\sum_{\L\cap\D(\f)\not=\emp}\F(\b,\L)
  -\sum_{\st\L\cap\D(\f)\not=\emp\atop\st\L\sset V\inda}\F(\b,\L)
  =\sum_{\st\L\cap\D(\f)\not=\emp\atop\st\L\cap (V\inda)^c\not=\emp}\F(\b,\L).
$$
Thus if $\f\in B$, $\L\cap\D(\f)\not=\emp$ and $\L\cap (V\inda)^c\not=\emp$, necessarily 
$\diam{\L}\geq (L/2)[(1/2)-\a]$.
By \equ(100) we obtain
$$
 \eqalign{
  &\quad |W_L^{\infty}(\b,\f)-W_L^{L/2}(\b,\f)|
   \leq\sum_{\st\L\cap\D(\f)\not=\emp\atop\st\L\cap (V\inda)^c\not=\emp}|\F(\b,\L)|\leq\cr
  &\leq\sum_{p\in V_L^{L/2}}\sum_{\st\L\ni p\atop\st\diam{\L}\geq(L/2)(1/2-\a)}|\F(\b,\L)|
  \leq L^2\nep{-m(\b)(L/2)(1/2-\a)}\leq C_1(\b).\cr
 }
$$
This bound and \equ(v8.2) give
$$
{\m\inda(\f|B)\over\bar\m\inda(\f)}
  \leq\nep{C_1}\bar\m\inda(B)^{-1}.
$$
To complete the proof we have to bound $\bar\m\inda(B)$ from below.
We refer to [2] to prove that there exists $C_2(\b)>0$ such that 
$\bar\m\inda(B)\geq \a C_2(\b)$.
\qed
\fine


\expandafter
 \ifx\csname sezioniseparate\endcsname\relax
  \fi

\numsec=-2
\numfor=1
\numtheo=1
\pgn=1

\beginsection B. References

\frenchspacing

\item{[1]}
     P. Diaconis and L. Saloff-Coste:
     Comparison theorems for reversible Markov chains.
     {\it Ann. Appl. Probab.} {\bf 3} (1993), no. 3, 696--730.
\item{[2]} 
   R. Dobrushin, R. Koteck\'y and S.~Shlosman:
   {\it Wulff Construction. A Global Shape from Local Interaction},
   Translation of Mathematical Monographs, {\bf 104} (1992).
   AMS.
\item{[3]}
   G. F. Lawler and A. D. Sokal:
   Bounds on the $L^2$ Spectrum for Markov Chains and Markov Processes: a Generalization of Cheeger's Inequality.
   {\it Trans. Amer. Math. Soc.} {\bf 309} (1988), no. 2, 557--580. 
\item{[4]} 
   T.~M.~Liggett:
   {\it Interacting particle systems}.
   Grundlehren der Mathematischen Wissenschaften {\bf 276} (1985).
   Springer-Verlag.
\item{[5]}    
  S.Lu and H.-T.Yau: 
  Spectral gap and logarithmic Sobolev inequality for Kawasaki and Glauber dynamics.
  {\it Comm. Math. Phys.} {\bf 156} (1993), no. 2, 399--433.
\item{[6]} 
   F.~Martinelli:
   On the two dimensional dynamical Ising model in the phase coexistence region.  
   {\it J. Statist. Phys.} {\bf 76} (1994), no. 5-6, 1179--1246.
\item{[7]} 
   F.~Martinelli:
   Lectures on Glauber dynamics for discrete spin models
   in {\it Lectures on probability theory and statistics (Saint-Flour, 1997)}.  
   Lecture Notes in Math. {\bf 1717} (1999), 93--191.
   Sringer--Verlag.
\item{[8]}
   F. Martinelli and E. Olivieri: 
   Approach to equilibrium of Glauber dynamics in the one phase region I: the attractive case.
   {\it Comm. Math. Phys.} {\bf 161} (1994), no. 3, 447--486.
\item{[9]}
   F. Martinelli and E. Olivieri:
   Approach to equilibrium of Glauber dynamics in the one phase region II: the general case.
   {\it Comm. Math. Phys.} {\bf 161} (1994), no. 3, 487--514.   
\item{[10]}
   G. Posta:
   Spectral Gap for an Unrestricted Kawasaki Type Dynamics.
   {\it ESAIM Probability \& Statistics} {\bf 1} (1997), 145--181.
\item{[11]}
   A.~D.~Sokal and L.~E.~Thomas:
   Absence of mass gap for a class of stochastic contour models.
   {\it J. Statist. Phys.} {\bf 51} (1988), no. 5-6, 907--947.
\item{[12]}
   D. W. Stroock and B. Zegarlinski: 
   The logarithmic Sobolev inequality for discrete spin on a lattice. 
   {\it Comm. Math. Phys.} {\bf 149} (1992), no. 1, 175--193.

\fine
\end